\documentclass[12pt, reqno]{amsart}

\usepackage[english]{babel}
\usepackage[latin1]{inputenc}
\usepackage{babel}
\usepackage{fancyhdr}
\usepackage{amsmath,amsfonts,amssymb,amsthm}
\usepackage[mathcal]{euscript}
\usepackage{mathrsfs}
\usepackage[tmargin=1.5in,bmargin=1.5in,lmargin=1.3in,rmargin=1.3in]{geometry}
\usepackage[all]{xy}
\usepackage{textcomp}
\usepackage{epsfig}
\usepackage{graphics}
\usepackage{graphicx}
\usepackage{mathrsfs}
\usepackage{caption}
\usepackage[OT2,OT1]{fontenc}
\newtheorem{thm}{Theorem}[section]

\newtheorem{lem}[thm]{Lemma}
\newtheorem{prop}[thm]{Proposition}
\newtheorem{defi}[thm]{Definition}

\newtheorem*{indhyp}{Induction hypothesis}

\usepackage{microtype}
\usepackage{bbm}
\usepackage{color}
\usepackage[textsize=small]{todonotes}       % includes TO DO LIST

\definecolor{halfgray}{gray}{0.55}
\definecolor{webgreen}{rgb}{0,.5,0}
\definecolor{webbrown}{rgb}{.6,0,0}
\definecolor{Maroon}{cmyk}{0, 0.87, 0.68, 0.32}
\definecolor{royalblue}{cmyk}{1, 0.50, 0, 0}
\definecolor{Black}{cmyk}{0, 0, 0, 0}

\usepackage{hyperref}
\hypersetup{%
    colorlinks=true, linktocpage=true, pdfstartpage=1, pdfstartview=FitV,%
    breaklinks=true, pdfpagemode=UseNone, pageanchor=true, pdfpagemode=UseOutlines,%
    plainpages=false, bookmarksnumbered, bookmarksopen=true, bookmarksopenlevel=1,%
    hypertexnames=true, pdfhighlight=/O,%nesting=true,%frenchlinks,%
    urlcolor=webbrown, linkcolor=royalblue, citecolor=webgreen, %pagecolor=RoyalBlue,%
}
\numberwithin{equation}{section}

\newcommand{\sss}{\scriptscriptstyle}

\renewcommand{\P}{\mathbb{P}}
\newcommand{\R}{\mathbb{R}}
\newcommand{\Pl}{\P_\lambda}

\newcommand{\indi}{\mathbbm{1}}
\newcommand{\Ccal}{\mathcal{C}}
\newcommand{\Cmax}{\Ccal_{\sss \mathrm{max}}}

\newcommand{\E}{\mathbb{E}}
\newcommand{\e}{\mathrm{e}}
\newcommand{\El}{\mathbb{E}_{\lambda}}

\begin{document}

\title[Connectivity Threshold for random  subgraphs of the Hamming graph]{Connectivity Threshold for random  subgraphs of the Hamming graph}

\author{Lorenzo Federico}

\author{Remco van der Hofstad}

\author{Tim Hulshof}
\email{l.federico@tue.nl, r.w.v.d.hofstad@tue.nl, w.j.t.hulshof@tue.nl}
\address{Department of Mathematics and Computer Science, Eindhoven University of Technology, PO Box 513, 5600 MB Eindhoven, The Netherlands}
\begin{abstract}
	We study the connectivity of random subgraphs of the $d$-dimensional Hamming graph $H(d, n)$, which is the Cartesian product of $d$ complete graphs on $n$ vertices. We sample the random subgraph with an i.i.d.\ Bernoulli bond percolation on $H(d,n)$ with parameter $p$. We identify the window of the transition: when $ np- \log n \to - \infty$ the probability that the graph is connected goes to $0$, while when $ np- \log n \to + \infty$ it converges to $1$.
	We also investigate the connectivity probability inside the critical window, namely when $ np- \log n \to t \in \mathbb{R}$. 
	We find that the threshold does not depend on $d$, unlike the phase transition of the giant connected component of the Hamming graph (see \cite{BorChaHofSlaSpe05a}). Within the critical window, the connectivity probability does depend on $d$. We determine how.\\
	
	\noindent{\sc Keywords:} \textit{connectivity threshold, percolation, random graph, critical window.}\\
\textit{MSC 2010:} 05C40, 60K35, 82B43. 
\end{abstract}

\date{\today}
\maketitle

\bigskip

\section{Introduction}

In this paper we  investigate the random edge subgraph of $d-$dimensional Hamming graphs. Hamming graphs are defined as follows:

\begin{defi}[Hamming graph] For integer $n$ write $[n] := \{1,\dots,n\}$. We define the $d-$dimensional Hamming graph $H(d,n)$ as the graph with vertex set
\[
V= [n]^d,
\]
and edge set
\[
E= \{ (v,w) : v,w \in V, \ v_j \neq w_j \text{ for exactly one }j \}. 
\]
\end{defi}
We study a percolation model on the Hamming graph. We define the random subgraph $H_{\lambda}(d,n)$ as the random edge subgraph with uniform edge retention probability $p=\frac{\lambda}{d(n-1)}$. Since the degree of every vertex in $H(d,n)$ is $d(n-1)$, the parameter~$\lambda$ thus indicates the expected number of outgoing edges from any given vertex. 

The phase transition for the existence of a giant component (i.e., when $|\Cmax| \approx \zeta |V|$ for $\zeta \in (0,1)$) was studied in  \cite{BorChaHofSlaSpe05a,HofNac12} for a larger class of finite transitive graphs that includes $H(d,n)$, while the slightly supercritical behavior was analyzed in \cite{HofLuc10} and \cite{HofLucSpe10} for $d=2$. 

In this work, we move away from the giant component critical point and we aim to determine the asymptotic probability that $H_{\lambda}(d,n)$ is connected for $d$ fixed and $n \to \infty$. The analogous problem was first solved for the \emph{Erd\H{o}s-R\'enyi Random Graph} (ERRG) in \cite{ErdRen59}. Observe that the ERRG arises as a special case of our problem if we put $d=1$. We will follow the proof for the ERRG (see e.g. \cite[Section 5.3]{Hofs09}), but we find that at places the internal geometry of the Hamming graph plays an important role. To overcome this difficulty we use an induction on the dimension $d$ and an exploration of the graph.

\section{Main Results}

Let $H_{n} :=H_{\lambda}(d,n)$ be a sequence of random edge subgraphs of $H(d,n)$ with parameter $\lambda=\lambda (n)$. Given $\lambda$ we want to determine the asymptotic probability that $H_n$ is connected.

\begin{thm}[Connectivity threshold for $H_\lambda (d,n)$]\label{main}

If $\lim_{n \to \infty} \lambda - d \log n=t \in \mathbb{R}$,
then 
\begin{equation} \label{crup}
\Pl(H_{n} \text {\emph{ is connected}}) \to \mathrm{e}^{-\mathrm{e}^{-t}}.
\end{equation}
Consequently, 
\begin{equation}
	\Pl(H_{n} \text {\emph{ is connected}}) \to  \begin{cases} 0 & \text{ if } \lambda - d \log n \to - \infty,\\
	1 & \text{ if } \lambda - d \log n \to + \infty.
	\end{cases}
\end{equation}

\end{thm}

These results show an interesting difference between the critical values of the giant component threshold and the connectivity threshold. The critical probability of the former, $p_{\sss GC} = \frac{1}{d(n-1)}(1+o(1))$, depends on $d$, while the latter, $p_{\sss \text{conn}} = \frac{\log n}{n-1}$, does not. This fact provides us with some insight into the structure of $H_\lambda(d,n)$ at the connectivity threshold: Consider the lower-dimensional ``hyperplanes'' (i.e., the subgraphs of $H(d,n)$ induced by all vertices $(v_1,\dots,v_d)$ that satisfy a set of constraints of the form $v_{j} = k_j$ for some $j \in [d]$, $k_j \in [n] $, see Definition \ref{def} below). Note that these hyperplanes are isomorphic to Hamming graphs of lower dimension. From \cite{BorChaHofSlaSpe05a} we know that there exist values of $\lambda$ such that $H_\lambda (d,n)$ has a giant component while the intersections of $H_\lambda  (d,n)$ with a hyperplane are subcritical (i.e., the largest components inside a hyperplane are of order $O(\log n)$). But an analogous property does not hold for the connectivity threshold: if $H_\lambda (d,n)$ is connected with probability converging to 1, then the same holds for all its hyperplanar subgraphs. 

We believe this phenomenon holds in much greater generality than Hamming graphs: our proof of Theorem \ref{main} can easily be modified to show that it also holds for the Cartesian product of $d$ copies of the complete $k$-partite graph, and we believe it to be true for a larger class of powers of high-degree transitive graphs.

\subsection{Related literature}
In \cite{ErdSpe79} Erd\H{o}s and Spencer studied the connectivity threshold of the \emph{hypercube} $H(d,2)$, where they found that the connectivity threshold occurs around $p=\frac{1}{2}$ (also independent of $d$). Clark \cite{Cla02} studied the connectivity threshold of $H(d,n)$ for $n$ fixed and $d \to \infty$, showing that if 
\[
p= 1 - \left( \dfrac{\xi(d)^{1/d}}{n} \right)^{\frac{1}{n-1}}
\]
and $\xi (d) \xrightarrow{d \to \infty} a \in (0,\infty)$, then the probability that the percolated graph is connected converges to $\e^{-a}$. Expansion of the above equation around $n = \infty$ shows that the $d \to \infty$ limit for large values of $n$ has the same behavior as the $n \to \infty$ limit.
Moreover, \cite{Joos13} shows that more generally, Cartesian products of fixed graphs have a connectivity threshold that only depends on their degree distribution as $d \to \infty$.

Sivakoff gives a statement analogous to our main theorem for site percolation in~\cite{Siv10}. It should be noted that site and edge percolation are very different models on the Hamming graph, as can be seen for instance in the fact that connectivity of site percolation on $K_n$ is trivial, whereas connectivity of edge percolation on $K_n$ (i.e., the ERRG) is not. See also \cite{Siva14}.

\section{Poisson convergence of isolated vertices}\label{sec3}

We start investigating the number of isolated vertices in the Hamming graph. As in the case of the ERRG, this provides a sharp lower bound on the window of the connectivity threshold. We define the number of isolated vertices
\begin{displaymath}
Y:= \sum_{i \in V} \indi_{\{|\Ccal_i|=1\}},
\end{displaymath}
where $\Ccal_i$ is the connected component of vertex $i$. We prove that in the critical window (i.e., when $\lambda - d \log n \to t \in \mathbb{R}$) the random variable $Y$
converges in distribution to a Poisson random variable. This proof is standard, and uses the same arguments applied to the proof given for the ERRG in \cite[Section 5.3]{Hofs09}.

Let $(x)_n$ denote the $n$th \emph{lower factorial} of $x$, i.e., $(x)_n := x (x-1) (x-2) \dotsm (x-n+1)$.
We will use the following lemmas, whose proofs are given in \cite[Section 2.1]{Hofs09} (for general versions see \cite[Chapter 6]{JanLucRuc00}):

\begin{lem}A sequence of integer-valued random variables $(X_n )_{n=1}^\infty$ converges in distribution to a Poisson random variable
with parameter $\mu$ when, for all $r = 1, 2, \dots,$

\begin{equation}\label{poi}
\lim_{n\to \infty} \E[(X_n)_r] = \mu^r.
\end{equation}
\end{lem}

\begin{lem}\label{lem:indi} When $X=\sum_{i \in \mathcal{I}}\indi_i$ is a sum of at least $r$ indicators then
\[
\E[(X)_r]=\sum_{i_1 \neq i_2 \neq ...\neq i_r} \P(\indi_{i_1}=\indi_{i_2}=\cdots=\indi_{i_r}=1),
\]
where the sum is over all sets of $r$ distinct indices.
\end{lem}
Given $H_n =(V_n,E_n)$, we want to prove that \eqref{poi} holds for $Y_n :=\sum_{v_i \in V_n} \indi_{\{|\Ccal_{v_i}|=1\}}$. We will use Lemma \ref{lem:indi} with an upper and lower bound on $\Pl(\indi_{i_1}=\dots =\indi_{i_r}=1)$ where we take $\indi_i$ to be the indicator function of the event that the vertex $v_i$ is isolated. Observe that we have $n^d !/(n^d - r)!$ different sets of distinct vertices of cardinality $r$. We call $m := d(n-1)$ the degree of $H(d,n)$.

The lowest probability comes from sets where none of the $r$ vertices are adjacent, hence we bound
\[
\Pl(\indi_{i_1}=\indi_{i_2}=\dots=\indi_{i_r}=1)\geq \left( 1- \dfrac{\lambda}{m}\right) ^{rm},
\]
while the highest probability comes from sets where all the $r$ vertices belong to the same $1$-dimensional subgraph, hence
\[
\Pl(\indi_{i_1}=\indi_{i_2}=\dots=\indi_{i_r}=1)\leq\left( 1- \dfrac{\lambda}{m}\right) ^{rm-\frac{r(r-1)}{2}}.
\]
For $n \leq r$ we can find better bounds but we do not mind, since we are interested in the asymptotic behavior when $n \to \infty$ and $r$ is fixed.
By the transitivity of the Hamming graph we bound, using $\lambda =d \log n + t(1+o(1))$,
\[\begin{split}
\El[(Y_n)_r] &\geq \dfrac{n^d !}{(n^d - r)!}\left( 1- \dfrac{\lambda}{m}\right) ^{rm}\\ &= \dfrac{n^d !}{(n^d - r)!} \e^{ -dr \log n-tr(1+o(1))}.
\end{split} \]
Since $\dfrac{n^d !}{(n^d - r)!}= n^{dr} (1-o(1))$, we find
\[
\El[(Y_n)_r] \geq n^{dr} \e^{ -dr \log n-tr} (1-o(1))= \e^{-tr(1+o(1))}.
\]
Similarly
\[\begin{split}
\El[(Y_n)_r] & \leq \dfrac{n^d !}{(n^d - r)!}\left( 1- \dfrac{\lambda}{m}\right) ^{rm-\frac{r(r-1)}{2}}\\ &= \dfrac{n^d !}{(n^d - r)!}n^{-dr} \e^{-tr(1 + o(1))} \left( 1- \dfrac{\lambda}{m}\right) ^{-\frac{r(r-1)}{2}}\\
& =  \e^{-tr(1 + o(1))}.
\end{split} \]

This proves that for each $r$, $\E[(Y_n)_r] \to  \e^{-tr}$ so that  by Lemma~\ref{poi} the distribution of $Y_n$ converges to $\text{Poi}(\e^{-t})$ when $\lambda - d \log n \to t$, so that

\begin{equation}\label{critical}
\Pl (Y_n=0) \to \e^{-\e^{-t}}.
\end{equation}

Furthermore $\{H_n \text{ connected}\} \subseteq  \{ Y_n=0 \}$, so we conclude that for $\lambda -d \log n \to t$
\[
\limsup_{n \to \infty} \Pl (H_\lambda (d,n) \text{ is connected}) \leq \e^{-\e^{-t}}.
\]

It remains to prove the matching lower bound, i.e., that in the critical window
\[
	\Pl(H_\lambda (d,n) \text{ is disconnected}\mid Y_n=0) \to 0.
\]

\section{Connectivity conditioned on no isolated vertices}\label{sec4}

We prove \eqref{crup} via induction on $d$. (The standard ``tree counting'' proof for the ERRG given in \cite[Section 5.3]{Hofs09} is too involved in the presence of geometry.)

\begin{indhyp}
If $\lim_{n \to \infty} \lambda - (d-1) \log n = t \in \R,$ then
\[
	\Pl(H_\lambda(d-1,n) \text{\emph{ is connected}}) \to \e^{-\e^{-t}},
\]
i.e., \eqref{crup} holds for $H (d-1,n)$.
\end{indhyp}

We initialize the induction by noting that $H(1,n)$ is a complete graph, so the random subgraph $H_\lambda (1,n)$ has the same distribution as an ERRG with $p=\frac{\lambda}{n-1}$. For this case it is proved in \cite{ErdRen59} that \eqref{crup} holds.

\begin{defi}[Hyperplanes]\label{def}
Given $H(d,n)=(V,E),$ define the \emph{hyperplanes}  $G_{jk}=(V_{jk},E_{jk})$ for some $j \in [d]$ and $k \in [n]$ as
\begin{itemize}
\item $V_{jk}= \{ (i_1 , i_2,\dots,  i_d) \in V : i_j=k \}$;
\item $E_{jk}= \{(v,w) \in E : v,w \in V_{jk}\}$.
\end{itemize}
Note that $H(d,n)$ has exactly $dn$ hyperplanes and that they are all isomorphic to $H(d-1,n)$.

We define $G_{jk}^\lambda$ as the intersection of the Random Edge Subgraph $H_\lambda (d,n)$ with the hyperplane $G_{jk}$, for each pair $j,k$.\end{defi}

The crucial idea of our proof is to show that once we have enough internally connected hyperplanes, all the remaining non-isolated vertices are connected to these connected hyperplanes with high probability. 

To use this argument, we condition on the event that a certain set of hyperplanes is internally connected. To ensure independence under this conditioning, we use disjoint edge sets to create the connected hyperplanes and to connect the remaining non-isolated vertices to them.

We define the sets $L= \{1,2,\dots, \lfloor n/2 \rfloor \} $ and $R = \{\lfloor n/2\rfloor +1, \lfloor n/2\rfloor+2,\dots,n \}$. For each $j \in [d]$, we divide $V$ into two sets:
\[
V_L (j) := \{ v \in V : v_{j} \in L\} \qquad \text{ and } \qquad V_R (j) := \{ v \in V : v_{j} \in R\}.
\]
This induces a partition on the edge set $E$:
\begin{equation*}
	\begin{split} E_{LL} (j) & := \{ (v,w) \in E : v,w \in V_L (j)\},\\
	 E_{RR} (j) & := \{ (v,w) \in E : v,w \in V_R (j)\},\\
E_{LR} (j) & := \{ (v,w) \in E : v \in V_L (j), w \in V_R (j)\}.
	\end{split}
\end{equation*}
For each fixed $j$ these sets are disjoint, so the occupation status of  the edges in one set is independent from the occupation status of edges in the other two sets.
Note that due to the geometry of $H (d,n)$ the exact composition of the sets $L$ and $R$ is not relevant, only their size matters. 

For some fixed $\alpha$ (to be determined later) we define the events
\begin{eqnarray*}
B_{R}(j) := \{G_{jk}^\lambda \text{ is connected for more than } \tfrac12 \alpha n \text{ different } k \in R \}, \\
B_L (j) := \{G_{jk}^\lambda \text{ is connected for more than } \tfrac12 \alpha n \text{ different } k \in L \}.
\end{eqnarray*}
We define $B := \bigcap_{j \in [d]} (B_L (j) \cap B_R (j))$.
In the final steps of the proof, on page \pageref{bsure}, we will show that $\Pl(B) \to 1$.

Note that the event $B$ states that there exist non-parallel internally connected hyperplanes, so when $B$ occurs, the geometry of the Hamming graph then ensures that \emph{all} internally connected hyperplanes are in the same connected component, deterministically. We will not make explicit use of this fact.
%By the inductive hypotesis and Weak Law of Large Numbers, if we choose $\alpha =\e^{-\e^{-(d-1)t/d}}- \varepsilon$, then for some $\varepsilon > 0$
%\begin{equation}
%\P (B_L (j)) \to 1, \quad \P (B_R (j)) \to 1.
%\end{equation}
Instead, we now prove that on $B$, with high probability, in the critical window $H_\lambda (d,n)$ consists only of the giant component and isolated points.

\begin{prop}\label{alfacrit} 
Let $\lambda - d\log n \to t \in \R$ and  $d \geq 2$, and let $\mathcal{I}$ be the set of all isolated points. Then
\[
\lim_{n \to \infty} \Pl(\{(\Cmax \cup \mathcal{I}) \neq [n]^d\} \cap B)= 0.
\]
\begin{proof}
We have to prove that with probability converging to 1 all edges present in the graph are connected to the giant component. We know that $|E| = \tfrac12 d n^d (n-1).$ We write $Z$ for the number of edges that do not connect to the giant component. If $Z=0$, then the claim holds, since all points outside the giant component must be isolated. We will prove that indeed $\El [Z \mid B] \to 0$.

Choose an edge $(v,v') \in E$ and let $i \in [d]$ be the unique direction such that $v_i \neq v'_i$. Choose $j \in [d]$ with $j \neq i$ and apply the partition defined above. Suppose that $v_j = v'_j \in L$ (the argument for $v_j = v'_j \in R$ is identical).
Define the event
\[
	F:= \{ v,v' \text{  are not connected to any internally connected hyperplane} \}.
\]
Since $B \subset B_R (j)$ for all $j$ it follows that 
\[
	\Pl(F \cap B) \le  \Pl(F \cap B_R (j)) \leq \Pl(F \mid B_R (j)).
\]

To bound $\Pl(F \mid B_R (j))$ we explore the graph starting from the vertices $v, v'$, with the following algorithm:
\begin{description}
\item[\texttt{Step 1}] Given the edge $(v , v')$, set as active the two end vertices: $a = v$, $a' = v'$. Initialize the set of searched vertices $\mathcal{S}= \varnothing$.
\item[\texttt{Step 2}] Set $\mathcal{S}= \mathcal{S} \cup \{a,a'\}$. Check all the edges $(a_1 , w)$, $(a_2 , w)$ such that $w$ belongs to a connected hyperplane $G^\lambda_{jr}$ with $r \in R$. 
\texttt{If} they are all vacant, \texttt{then} go to Step 3, \texttt{else} stop.
\item[\texttt{Step 3}] Check all the edges $(a , w)$, $(a' , w)$ such that $w \in V_L (j)\setminus \mathcal{S}$. We define the set $W=\{ w \in V_L (j)\setminus \mathcal{S}: $ one or both of $(a , w)$ and $(a' , w)$ are occupied$\}$.
\texttt{If} $|W| \geq 2$ \texttt{then} go to Step 4, \texttt{else} stop.
\item[\texttt{Step 4}] Choose $w$, $w' \in W$ according to an arbitrary but fixed rule and set them as the active vertices: $a=w$ and $a'= w'$. Return to Step 2.
\end{description}

Activating only two vertices at each cycle of the algorithm allows for some control over the depletion of points outside the connected hyperplanes. This means that the algorithm can terminate before the starting edge has been connected to the giant component or before its connected component has been completely explored. This is not a problem: the calculations below show that this algorithm gives a sufficiently sharp result to prove the claim.

Indeed, we want to show that the probability that the exploration process terminates before finding a connection to one of the internally connected hyperplanes among $\{G^\lambda_{jr} \, \colon \, r \in R\}$ has probability tending to zero. That is, we want to show that with high probability, the algorithm does not terminate during Step 3. We write $T$ for the cycle at which this happens. We set $T = \infty$ if the process finds the giant component, namely if the algorithm terminates during Step 2.

Note that the algorithm is designed with certain independencies. Indeed, the event $B_R (j)$ depends only on the edges in $E_{RR} (j)$, Step 2 of the exploration only depends on edges in $E_{LR} (j)$, and Step 3 only depends on edges in $E_{LL} (j)$.

Let $P_g=P_g(k)$ be the probability of finding a connection to one of the internally connected hyperplanes among $\{G^\lambda_{jr} \, \colon \, r \in R\}$ 
%the giant component in $R$ 
during the $k$-th cycle of the exploration algorithm, conditioned on the event $B_R (j)$ and on the event that the algorithm has not yet terminated. We bound
\[
1- P_g \leq \left( 1- \frac{\lambda}{m}\right)^{\alpha n} = \e^{-\lambda \alpha/d}(1+o(1)) \le C n^{-\alpha},
\]
for a constant $C$ that depends on $t$.
(This bound does not depend on $k$ because the algorithm terminates as soon as the exploration finds a connected hyperplane, so there is no depletion of points inside the connected hyperplanes.)

Let $N_k$ denote the number of vertices discovered in Step 3 of the $k$-th cycle of the exploration and let 
\[
	P_{k,2} := \Pl(N_k \geq 2 )=1 - \Pl (N_k =0) - \Pl (N_k =1).
\]

Each vertex $v \in V_L (j)$ has $(d-1)(n-1)+\lfloor n/2\rfloor$ neighbors in $V_L (j)$, and at time~$k$ at most $2k$ of them have already been explored, so $N_k$ stochastically dominates a Bin$ ((2d-1)(n-1) -4k, \lambda /m)$ random variable. We bound
\[
\Pl(N_k \in \{ 0,1\} )\leq  \left(1-\frac{\lambda}{m}\right)^{(2d-1 ) (n-1)-4k} +(2d-1 ) n\frac{\lambda}{m} \left(1-\frac{\lambda}{m}\right)^{(2d-1 ) (n-1)-4k -1}.
\]
So we obtain for some constant $c$
\[\begin{split}
1-P_{k,2} \leq  & \left(1-\frac{\lambda}{m}\right)^{(2d-1 ) (n-1)-4k}\left(1+\dfrac{2d-1}{d} \frac{\lambda}{1-\frac{\lambda}{m}}\right) \\
\leq & c \lambda n^{-2d+1}\left(1-\frac{\lambda}{m}\right)^{-4k}.
\end{split}\]

If $T=s$, then the exploration does not reach a connected hyperplane during the first $s$ cycles and then the algorithm terminates on Step 3 of the $s$-th cycle, so we can bound
\[\begin{split}
	\Pl(T = s\mid B_R (j))&\leq (1-P_{s,2}) \prod_{k \le s} (1-P_g) \\
			&\leq C \lambda n^{-2d+1}\left(1-\frac{\lambda}{m}\right)^{-4s} n^{-\alpha s}\\ 			&=  C \lambda n^{-2d+1}\left(\left(1-\frac{\lambda}{m}\right)^4 n^{\alpha}\right)^{-s}  ,
\end{split}\]
for some constant $C$ that depends on $t$.
It follows that
\[
\Pl(T < \infty \mid B_R (j)) \le  C \lambda n^{-2d+1 } \sum_{s=1}^\infty \left(\left(1-\frac{\lambda}{m}\right)^4 n^{\alpha}\right)^{-s}.
\]
For sufficiently large $n$ we have $n^{\alpha }(1-\frac{\lambda}{m})^4>1$, so the tail of the sum behaves like a convergent geometric series, and we bound
\[
	\Pl(T < \infty \mid B_R (j)) \le C n^{-2d+1-\alpha}\log n \ll  \frac{2}{d}n^{-d}(n-1)^{-1},
\]
for all $d \geq 2$.

Since $H(d,n)$ is transitive and there are  $\frac{d}{2}n^{d}(n-1)$ edges, 
\[
\El [Z \mid B] \leq \frac{d}{2}n^{d}(n-1) \Pl(T < \infty \mid B_R (j)) \to 0,
\]
and the claim now follows.
\end{proof}
\end{prop}

\subsection*{Completion of the proof: induction on the dimension}
Recall that the case $d=1$ initiates the induction, since $H_\lambda (1,n)$ is an Erd\H{o}s-R\'enyi graph, so \eqref{crup} holds.

For the inductive step we assume that \eqref{crup} holds for $H_\lambda (d-1,n)$, i.e., that for all $t \in \mathbb{R}$ and all sequences $\lambda = \lambda(n)$ such that $\lim_{n \to \infty} \lambda - (d-1) \log n = t$, we have
\[
	\Pl(H_\lambda(d-1,n) \text{ is connected}) \to  \e^{-\e^{-t}}.
\]
We want to prove that the same holds for $H_\lambda (d,n)$. 

Given $H_\lambda (d,n)$, its intersection $G_{jk}^\lambda$ with the hyperplane $G_{jk}$ has the same distribution as $H_{\frac{d-1}{d}\lambda}(d-1,n)$ since $p= \frac{\lambda}{d(n-1)}$, and each vertex has $(d-1)(n-1)$ outgoing edges in $G_{jk}$.
We assumed that $\lim_{n \to \infty} \lambda - d \log n = t$, which implies that
\[
	\lim_{n \to \infty} \frac{d-1}{d} \lambda - (d-1) \log n = \frac{d-1}{d} t.
\]

Note moreover that $\indi_{\{G_{jk}^\lambda \text{ is connected}\}}$ and $\indi_{\{G_{jk'}^\lambda \text{ is connected}\}}$ are i.i.d.\ random variables when $k \neq k'$ under $\Pl$, since for fixed $j$ all the subgraphs $G_{jk}^\lambda$ are i.i.d.\ random subgraphs with the same law as $H_{\frac{d-1}{d}\lambda}(d-1,n)$. It thus follows by the inductive hypothesis that the asymptotic probability that $G_{jk}^\lambda$ is connected is $\exp(-\e^{-(d-1)t/d})$.

If we choose $\varepsilon > 0$ such that $\alpha := \exp(-\e^{-(d-1)t/d})- \varepsilon > 0$, then for each $j$ the induction hypothesis and the Weak Law of Large Numbers imply that that
\[
\begin{split}
\Pl (B_L (j)^c) & = \Pl \Biggl( \sum_{k=1}^{\lfloor n/2 \rfloor} \indi_{\{G_{jk}^\lambda \text{ is connected}\}} \leq \tfrac12 \alpha n\Biggr)\\
& \leq  \Pl \Biggl( \Biggl|  \dfrac{2}{n}\sum_{k=1}^{\lfloor n/2 \rfloor} \indi_{\{G_{jk}^\lambda \text{ is connected}\}} - \e^{-\e^{-(d-1)t/d}}\Biggr|  > \varepsilon\Biggr) \to  0 \ \ \text{ as } n \to \infty .
\end{split}
\]
The same is true for each $B_R (j)$.
Using the union bound,
\begin{equation}\label{bsure}
\Pl (B^c) \leq \sum_{j=1}^d \Big(\Pl (B_L (j)^c)+ \Pl( B_R (j)^c)\Big) \to 0. 
\end{equation}

Finally, we combine \eqref{critical}, \eqref{bsure} and Proposition \ref{alfacrit} to obtain
\[\begin{split}
\Pl (Y>0) & \leq  \Pl (H_\lambda(d,n) \text{ is disconnected}) \\ 
& \leq  \Pl (Y>0) + \Pl (B^c) +\Pl(\{(\Cmax \cup \mathcal{I}) \neq [n]^d\} \cap B),
\end{split}\]
completing the proof of the main theorem. \qed

\section*{Acknowledgments}
We would like to thank the anonymous referee for the very helpful comments.

The work of LF, RvdH and TH is supported by the Netherlands Organisation for Scientific Research (NWO) through the Gravitation {\sc Networks} grant 024.002.003. The work of RvdH is also supported by the Netherlands Organisation for Scientific Research (NWO) through VICI grant 639.033.806.

%%%%%%%%%%%%%%%%%% REFERENCES %%%%%%%%%%%%%%%%%%%%%%%%%%%%%%%%%%%%%%%%%%%%%%%%%%%
\begin{small}
\bibliographystyle{abbrv}
\bibliography{LorenzosBib}

\end{small}

\end{document}